\documentclass[oneside,leqno,12pt]{article}
\usepackage{amssymb,amsmath,latexsym}

\def\C{\mathbb{C}}

\def\N{\mathbb{N}}

\def\P{\mathbb{P}}
\def\Q{\mathbb{Q}}
\def\R{\mathbb{R}}

\def\Z{\mathbb{Z}}

\def\cC{\mathcal{C}}

\def\cF{\mathcal{F}}
\def\cG{\mathcal{G}}
\def\cH{\mathcal{H}}

\def\cL{\mathcal{L}}
\def\cFl{\mathcal{F}\ell}
\def\fl{\it{fl}}

\def\cO{\mathcal{O}}

\def\inj{\varinjlim}
\def\proj{\varprojlim}

\def\id{{\rm id}}

\def\nek{\text{\hbox{$\simeq$ \kern-.95em \hbox{$/$ \kern.05em}}}}

\DeclareMathOperator{\codim}{codim}

\DeclareMathOperator{\Hom}{Hom}

\DeclareMathOperator{\Span}{span}

\DeclareMathOperator{\Pic}{Pic}

\def\cplus{\hbox{$\subset${\raise1.05pt\hbox{\kern -0.55em
${\scriptscriptstyle +}$}}\ }}
\def\bcplus{\hbox{$\supset${\raise1.05pt\hbox{\kern -0.55em
${\scriptscriptstyle +}$}}\ }}

\def\ctimes{\hbox{$\times${\raise1.1pt\hbox{\kern -0.27em
${\scriptscriptstyle |}$}}\ }}

\def\udarrow{\hbox{$\nearrow${\kern -0.97em$\searrow$}\ }}

\def\bctimes{\hbox{$\times${\raise1.1pt\hbox{\kern -.74em
${\scriptscriptstyle |}$}}\ }\,\,}

\newtheorem{theorem}{Theorem}
\newtheorem{proposition}{Proposition}

\title  {Ind--varieties of generalized flags as homogeneous spaces
for classical ind--groups}
\author {Ivan Dimitrov and Ivan Penkov}
\date{\empty}

\begin{document}

\maketitle

\begin{abstract}
The purpose of the present paper is twofold: to introduce the notion of a generalized flag
in an infinite dimensional vector space $V$ (extending the notion of a flag of subspaces in a 
vector space), and to give a geometric realization of homogeneous spaces
of the ind--groups $SL(\infty)$, $SO(\infty)$ and $Sp(\infty)$ in terms of generalized flags.
Generalized flags in $V$ are chains of subspaces which 
in general cannot be enumerated by integers.
Given a basis $E$ of $V$, we define a notion of $E$--commensurability for generalized flags, 
and prove that the set
$\cFl (\cF, E)$ of generalized flags $E$--commensurable with a fixed generalized flag $\cF$ in $V$ has 
a natural structure of an ind--variety. In the case when $V$ is the standard representation
of $G = SL(\infty)$, all homogeneous ind--spaces $G/P$ for parabolic subgroups $P$
containing a fixed splitting Cartan subgroup of $G$, are of the form $\cFl (\cF, E)$.
We also consider isotropic generalized flags. The corresponding ind--spaces are homogeneous spaces for 
$SO(\infty)$ and $Sp(\infty)$. As an application of the construction,
we compute the Picard group of $\cFl (\cF, E)$ (and of its isotropic analogs)
and show that $\cFl (\cF, E)$ is a projective ind--variety if
and only if $\cF$ is a usual, possibly infinite, flag of subspaces in $V$.

Key words (2000 MSC): Primary 14M17, Secondary 14L35, 14M15.
\end{abstract}

\vskip-6in 
\hskip3.75in
{\small To Joseph A. Wolf with admiration}
\vskip6in

\section*{Introduction} \label{Sec0}
Flag varieties play a fundamental role both in representation theory and algebraic geometry.
There are two standard approaches to flag varieties: the group--theoretic one, where a 
flag variety is defined as $G/P$ for a classical algebraic group $G$ and a parabolic
subgroup $P$, and the geometric one, where a flag variety is defined as the set of all chains of
subspaces of fixed dimensions in a finite dimensional vector space, which in addition are
assumed isotropic in the presence of a bilinear form. The very existence of these two approaches
is in the heart of the interplay between representation theory and geometry.

The main topic of this paper is a purely geometric construction of homogeneous spaces
for the classical ind--groups $SL(\infty)$, $SO(\infty)$ and $Sp(\infty)$.
Despite the fact that many 
phenomena related to inductive limits of classical groups have been studied
(see for instance,  \cite{O1}, \cite{O2}, \cite{SV}, \cite{VK1}, \cite{VK2},
many natural questions remain unanswered. In particular, 
the only approach to homogeneous spaces of classical ind--groups discussed in the literature is a
representation--theoretic one, and has been introduced by Joseph A. Wolf and his collaborators,
\cite{NRW}, \cite{DPW}. 

The difficulty in the purely geometric approach is that
the consideration of flags, i.e. chains of subspaces enumerated by integers, is no
longer sufficient. To illustrate the problem, 
let, more specifically, $G$ denote the ind--group $SL(\infty)$ over a field of characteristic $0$, and
$P$ be a parabolic subgroup of $G$. By definition, $G$ is the union
of a standard system of nested algebraic groups $SL(n)$ and $P$ is the
union of parabolic subgroups.
If $V$ is the natural representation of
$G$, all $P$--invariant subspaces in $V$ form a chain $\cC$
of subspaces of $V$. 
In general, the chain $\cC$ has a rather complicated structure and is not necessarily a flag,
i.e. cannot be indexed by integers.
We show, however, that $\cC$ always contains 
a canonical subchain $\cF$ of subspaces of $V$ with the property that
every element of $\cF$ is either the immediate predecessor $F'$ of a
subspace $F'' \in \cF$ or the immediate successor $F''$ of a subspace $F' \in \cF$, and, in addition,
each nonzero vector $v \in V$ belongs to a difference $F'' \backslash F'$. 
These two properties define generalized flags. (Maximal generalized flags already appeared 
in \cite{DP} in a related but somewhat different context.) 
If, in addition, the vector space $V$ is equipped with a non--degenerate bilinear
(symmetric or antisymmetric) form, we introduce the notion of an isotropic generalized flag.

Informally we think of two (possibly isotropic)
generalized flags being  commensurable if they only differ in a 
finite dimensional subspace of $V$ in which they reduce to flags of the same type.
The precise definition is given in Section 4.
The main result of this paper is the construction of the ind--varieties of commensurable 
generalized flags and their identification with homogeneous ind--spaces $G/P$ for
classical locally linear ind--groups $G$ isomorphic to $SL(\infty)$, $SO(\infty)$, or $Sp(\infty)$,
and corresponding parabolic subgroups $P$. 

The paper is completed by providing two applications: an explicit computation of the Picard group of
any ind--variety of commensurable generalized flags $X$ and a criterion for projectivity of $X$.
We show that the Picard group of $X$
admits a description very similar to the classical one;
however, $X$ is projective if and only if it is an ind--variety of usual flags. 

The ``flag realization'' of the ind--varieties $G/P$ given in the present paper
opens the way for a detailed and explicit study of the geometry of $G/P$, which should
play a role as prominent as the geometric representation theory of the classical
algebraic groups.

\noindent
{\bf Acknowledgement}
\newline
We thank Vera Serganova for a detailed and thoughtful
critique of the first version of this paper. The second named author thanks the Max Planck Institute
for Mathematics in Bonn for support and hospitality.

\noindent
{\bf Conventions}
\newline
$\N$ stands for $\{ 1, 2, \ldots\}$ and $\Z_+ = \N \cup \{0\}$.
The ground field is a field $k$ of characteristic $0$ which will be assumed
algebraically closed only when explicitly indicated in the text. As usual,
$k^\times$ is the multiplicative group of $k$.
The superscript $^*$ denotes dual vector space. The signs $\inj$ and $\proj$
stand respectively for direct and inverse limit over a direct or inverse system of morphisms 
parametrized by $\N$ or $\Z_+$.  $\Gamma(X, \cL)$ denotes the global sections of 
a sheaf $\cL$ on a topological space $X$.  All orders are assumed linear and
strict, and all partial orders are assumed to have the additional property that the relation
"neither $x \prec y$ nor $y \prec x$" is an equivalence relation.

\section{Preliminaries}
An {\em ind--variety} (over $k$) is a set $X$ with a filtration
\begin{equation} \label{ind_filt}
X_0 \subset X_1 \subset \ldots \subset X_n \subset X_{n+1} \subset \ldots
\end{equation}
such that
$X = \cup_{n\in \Z_+} X_n$, each $X_n$ is a Noetherian algebraic variety, and
the inclusions $X_n \subset X_{n+1}$ are closed immersions of algebraic varieties.
An ind--variety $X$ is automatically a topological space: a subset
$U \subset X$ is {\em open} in $X$ if and only if, for each $n$, $U \cap X_n$ is an
open subvariety of $X_n$\,.
The {\em sheaf of regular functions on} $X$, or {\it the structure sheaf}
$\cO_X$ of $X$, is the inverse limit
$\cO_X = \proj \cO_{X_n}$ of the sheaves of regular functions $\cO_{X_n}$ on the $X_n$\,.
An ind--variety $X = \cup_{n \in \Z_+} X_n = \inj X_n$ is {\em proper} if and only if all the
varieties $X_n$ are proper, is {\em affine} if and only if all the $X_n$ are affine.
A {\em morphism} from an ind--variety $X$ to an ind--variety $Y$ is a
map $\varphi : X \to Y$ such that, for every $n \geq 0$, the restriction
$\varphi|_{X_n}$ is a morphism of $X_n$ into $Y_m$ for some $m = m(n)$.
An {\em isomorphism} of ind--varieties is a morphism which admits an inverse
morphism.  An {\em ind--subvariety} $Z$ of $X$ is a subset $Z \subset X$
such that $Z \cap X_n$ is a subvariety of $X_n$ for each $n$.
Finally, an {\it ind--group} is by definition a group object in the category
of ind--varieties.
In this paper we consider only ind--groups $G$ which are {\it locally
linear}, i.e. ind--varieties $G$ with an ind--variety filtration $G_0 \subset G_1 \subset 
\ldots \subset G_n \subset G_{n+1} \subset \ldots$, such that all $G_n$ are linear algebraic 
groups and the inclusions are group morphisms.

Let $V$ be a vector space of countable dimension.  Fix an integer $l \geq 1$.
The set $Gr(l;V)$ of all $l$--dimensional subspaces of $V$ has a canonical structure
of proper ind--variety: any filtration
$0 \subset V_l \subset V_{l+1} \subset \ldots \subset V = \cup_{r \geq 0} V_{l+r}$,
$\dim V_{l+r} = l+r$, induces a filtration 
$$
Gr(l;V_l) \subset Gr(l;V_{l+1}) \subset \ldots \subset Gr(l;V),
$$
and the associated ind--variety structure on $Gr(l;V)$ is independent of
choice of filtration on $V$.  For $l = 1$, $\P(V) := Gr(1;V)$ is by definition the
{\em projective ind--space associated to} $V$.

An {\it invertible sheaf} on an ind--variety $X$ is a sheaf of $\cO_X$--modules locally
isomorphic to $\cO_X$. The set of isomorphism classes of invertible sheaves on $X$
is an abelian group (the group structure being induced by the operation of 
tensor product over $\cO_X$ of invertible sheaves). By definition, the latter
is the {\it Picard group $\Pic X$ of $X$}. It is an easy exercise to show that 
$\Pic X = \proj \Pic X_n$ for any filtration (\ref{ind_filt}). If $X = \P(V)$, then
$\Pic X \cong \Z$. The preimage of $1$ under this isomorphism is the class of the 
standard sheaf $\cO_{\P(V)}(1)$ where, by definition, 
$\cO_{\P(V)}(1) := \proj \cO_{\P(V_n)}(1)$.

An invertible sheaf $\cL$ on a proper ind--variety $X$ is {\it very ample} if, for some
filtration (\ref{ind_filt}), its restrictions $\cL_n$ on $X_n$ are very ample 
for all $n$, and all restriction maps $\Gamma (X_n; \cL_n) \to \Gamma (X_{n-1}, \cL_{n-1})$ 
are surjective. A very ample invertible sheaf defines a closed immersion of $X$ into 
$\P(\inj \Gamma (X_n, \cL_n)^*)$ as for each $n$ the restrictions $\cL_n$ and $\cL_{n-1}$
define a commutative diagram of closed immersions
$$
\begin{array} {ccc}
X_{n-1} & \hookrightarrow & \P(\Gamma (X_{n-1}, \cL_{n-1})^*)\\
\cap & & \cap \\
X_n & \hookrightarrow & \P(\Gamma (X_n, \cL_n)^*).
\end{array}
$$
Conversely, given a closed immersion $X \hookrightarrow \P(V)$, the inverse image
of $\cO_{\P(V)}(1)$ on $X$ is a very ample invertible sheaf on $X$. Therefore, a
proper ind--variety $X$ is projective, i.e. $X$ admits a closed immersion into 
a projective ind--space, if and only if it admits a very ample invertible sheaf.

\section{Generalized flags: definition and first properties} \label{Sec1}
Let $V$ be a vector space over $k$. A {\it chain of subspaces in $V$} is a set $\cC$ of pairwise distinct
subspaces of $V$ such that for any pair $F'$, $F'' \in \cC$, 
either $F' \subset F''$ or $F'' \subset F'$. Every chain of subspaces $\cC$
is ordered by proper inclusion. Given $\cC$, we denote by $\cC'$ (respectively, by $\cC''$) the subchain of $\cC$ which 
consists of all $C \in \cC$ with an immediate successor (respectively, an immediate predecessor).
A {\it generalized flag in $V$}
is a chain of subspaces $\cF$ which satisfies the following conditions:
\begin{itemize}
\item[]{(i)} each $F \in \cF$ has an immediate successor or an immediate predecessor, i.e. $\cF = \cF' \cup \cF''$;
\item[]{(ii)} $V \backslash \{0\} = \cup_{F' \in \cF'} F'' \backslash F'$, where $F'' \in \cF''$ is the immediate
successor of $F' \in \cF'$.
\end{itemize}
Given a generalized flag $\cF$ and a subspace $F'' \in \cF''$ (respectively, $F' \in \cF'$), we
will always denote by $F'$ (resp., by $F''$) its immediate predecessor (resp., immediate successor).
Furthermore, condition (ii) implies that each nonzero vector $v \in V$ determines a unique pair
$F_v' \subset F_v''$ of subspaces in $\cF$ with $v \in F_v'' \backslash F_v'$.

\noindent
{\bf Example 1.}

\noindent
(i) We define a {\it flag in $V$} to be a chain of subspaces $\cF$ satisfying (ii)
and which is isomorphic as an ordered set to a subset of $\Z$. 
A flag can be equivalently defined as a chain of subspaces $\cF$ for which there 
exists a strictly monotonic map of ordered sets $\varphi: \cF \to \Z$ and, in addition, 
$\cap_{F \in \cF} F = 0$ and $\cup_{F \in \cF} = V$. There are four different kinds of flags:
a finite flag of length $k$ $\cF = \{0 = F_1 \subset F_2 \subset \ldots \subset F_{k-1} \subset F_k = V\}$;
an infinite accending flag $\cF = \{0 = F_1 \subset F_2 \subset F_3 \ldots\}$, where $\cup_{i \geq 1} F_i = V$;
an infinite descending flag $\cF = \{ \ldots \subset F_{-3} \subset F_{-2} \subset F_{-1} = V \}$, where 
$\cap_{i \leq -1} =0$; and a two sided infinite flag $\cF = \{ \ldots \subset F_{-1} \subset F_0 \subset F_1 \subset \ldots\}$,
where $\cap_{i \in \Z} F_i = 0$ and $\cup_{i \in \Z} F_i = V$.  

\noindent
(ii) One of the simplest examples of a generalized flag in $V$ which is not a flag is a generalized flag with
both an infinite accending part and an infinite descending part, i.e. 
$\cF = \{0 = F_1 \subset F_2 \subset F_3 \subset \ldots \subset F_{-3} \subset F_{-2} \subset F_{-1} = V\}$,
where $\cup_{i \geq 1} F_i = \cap_{j \leq -1} F_j$. 

\noindent
(iii) Let $V$ be a countable dimensional vector space with basis
$\{e_q\}_{q \in \Q}$. Set $F'(q):= \Span\{e_r \, | \, r < q\}$ and 
$F''(q):= \Span\{e_r \, | \, r \leq q\}$. Then 
$\cF := \cup_{q \in \Q} \{F'(q), F''(q)\}$ is a generalized flag in $V$.
$\cF$ is not a flag, and moreover, no $F \in \cF$ has both an immediate predecessor
and an immediate successor.

The following Proposition shows that each of the subchains
$\cF'$ and $\cF''$ reconstructs $\cF$.

\begin{proposition} \label{pro01}
Let $\cF$ be a generalized flag in $V$. Then
\newline
{\rm (i)} for every $F' \in \cF'$, $F' = \cup_{G'' \in \cF'', G'' \subset F'', G'' \neq F''} G''$;
\newline
{\rm (ii)} for every $F'' \in \cF''$, $F'' = \cap_{G' \in \cF', G' \supset F', G' \neq F'} G'$.
\end{proposition}

\noindent
{\bf Proof.} (i) The inclusion 
$F' \supset \cup_{G'' \in \cF'', G'' \subset F'', G'' \neq F''} G''$ is obvious.
Assume now that $v \in F'$. Let $v \in H'' \backslash H'$ for some $H' \in \cF'$ and its immediate 
successor $H'' \in \cF''$.
Then $H' \subset F'$ and hence $H'' \subset F'$, i.e.
$v \in \cup_{G'' \in \cF'', G'' \subset F'', G'' \neq F''} G''$
which proves that $F' \subset \cup_{G'' \in \cF'', G'' \subset F'', G'' \neq F''} G''$.
Assertion (ii) is proved in a similar way.
\hfill $\square$

Any chain $\cC$ of subspaces in $V$ determines the 
following partition of $V$:
\begin{equation} \label{eq01}
V = \sqcup_{v \in V} [v]_\cC, \quad {\text { where }} \quad
[v]_\cC := \{ w \in V \, | \, w \in F \Leftrightarrow v \in F, \,\, \forall F \in \cC\}.
\end{equation}
Consider this correspondence as a map $\pi$ 
from the set of chains of subspaces in $V$ into the set
of partitions of $V$.
This map is not injective, for $\pi(\cC') = \pi(\cC)$ if $\cC'$ is obtained from $\cC$  
by adding arbitrary intersections and unions of elements of $\cC$.
As we show in Proposition \ref{pro02} below, 
the notion of a generalized flag provides us with a natural
right inverse of $\pi$, i.e. with a map $\gamma$ (defined on the image of $\pi$)
such that $\pi \circ \gamma = \id$. 
This explains the special role of generalized flags among arbitrary chains of subspaces in
$V$. Namely, every generalized flag in $V$ is a natural representative of 
the class of chains of subspaces in $V$ which yield the same partition of $V$.

\begin{proposition} \label{pro02}
Given a chain $\cC$ of subspaces in $V$, there exists a unique generalized flag $\cF$ in $V$
for which $\pi(\cC) = \pi(\cF)$.
\end{proposition}

\noindent
{\bf Proof.}
To prove the existence, set $F_v' := \cup_{W \in \cC, v \not \in W} W$ and 
$F_v'' := \cap_{W \in \cC, v \in W} W$, and put $\cF := \cup_{v \in V \backslash \{0\}} \{ F_v', F_v''\}$.
It is obvious from the definition of $\cF$ that $\pi(\cC) = \pi(\cF)$. To show that $\cF$ is a generalized flag,
notice that, for any pair of nonzero vectors $u, v \in V$,
exactly one of the following three possibilities holds:
\begin{itemize}
\item[] $F_u' = F_v'$, and hence $F_u'' = F_v''$;
\item[] $F_u' \subset F_v'$, and hence $F_u'' \subset F_v'$;
\item[] $F_u' \supset F_v'$, and hence $F_u' \supset F_v''$.
\end{itemize}
Indeed, if, for every $W \in \cC$, $u \in W$ if and only if $v \in W$, then 
$F_u' = F_v'$ and $F_u'' = F_v''$. Assume now, that there exists $W \in \cC$ such that
$u \in W$ but $v \not \in W$. Then $F_u'' \subset W \subset F_v'$. Similarly, if there
exists $W \in \cC$ such that
$u \not \in W$ but $v \in W$, we have $F_v'' \subset W \subset F_u'$.
The existence of
$\cF$ is now established.

The uniqueness follows from the fact that $[v]_{\cC} = (\cap_{W \in \cC, v \in W} W) \backslash
(\cup_{W \in \cC, v \not \in W} W)$, while, for a generalized flag $\cF$, $[v]_{\cF} = F_v'' \backslash F_v'$.
\hfill $\square$

We now define the map $\gamma$ by setting $\gamma (\pi (\cC)):= \cF$, and put $\fl := \gamma \circ \pi$.
In the example below we determine
the preimages under $\fl$ of the generalized flags introduced in Example 1.  
The computation is based on the following simple fact: if $\bar{\cC}$ is any chain in $\fl^{-1}(\cF)$, 
then every nonzero subspace $\bar{C} \in \bar{\cC}$ is the union of spaces from $\cF$.

\noindent
{\bf Example 2.} The cases (i), (ii) and (iii) below refer to the corresponding cases in Example 1.

\noindent
(i) If $\cF$ is a flag in $V$ then $\fl^{-1} (\cF)$ consists of $\cF$ and the chains obtained from $\cF$
by adding $0$, $V$ or both, in case $0$ and/or $V$ do not belong to $\cF$.

\noindent
(ii) In this case $\fl^{-1}(\cF)$ consists of two chains: $\cF$ itself and the chain obtained by adding
$\cup_{i \geq 1} F_i = \cap_{j \leq -1} F_j$ to $\cF$.

\noindent
(iii) In this case there are infinitely many chains $\bar{\cC}$ with $\fl(\bar{\cC}) = \cF$.
Set $F'(x) := \Span \{e_r \, | \, r < x\}$ for any $x \in \R$,
and let $\cC$ denote the chain $\{F'(x) \, | \, x \in \R\} 
\cup \{F''(q) \, | \, q \in \Q\} \cup \{ 0 , V\}$. It is easy to check that $\fl(\cC) = \cF$ and that
any chain in $\fl^{-1}(\cF)$ is a subchain of $\cC$. To characterize explicitly all chains in 
$\fl^{-1}(\cF)$, for any subchain $\bar{\cC} \subset \cC$, set
$\R_{\bar{\cC}} := \{ x \in \R \, | \, F'(x) \in \bar{\cC} \}$ and 
$\Q_{\bar{\cC}} := \{ q \in \Q \, | \, F''(q) \in \bar{\cC} \}$.
Then $\fl(\bar{\cC}) = \cF$ if and only if, for any $r \in \Q$, we have 
$r \in \Q_{\bar{\cC}}$ or $r = \inf \{x \in \R_{\bar{\cC}} \cup \Q_{\bar{\cC}} \, | \,   r < x \}$, and
$r \in \R_{\bar{\cC}}$ or $r = \sup \{x \in \R_{\bar{\cC}} \cup \Q_{\bar{\cC}} \, | \, x < r\}$.

A generalized flag $\cF$ in $V$ is {\it maximal} if it is not properly contained in another
generalized flag in $V$. It is easy to see that 
the generalized flags introduced in Example 1(ii) and (iii) are 
maximal. More generally, a generalized flag 
$\cF$ is maximal if and only if $\dim(F_v''/F_v') = 1$ for every nonzero $v \in V$. 
Indeed, assume $\dim(F_{v_0}''/F_{v_0}') > 1$ for some $v_0$. Let $F \subset V$ be a subspace
with proper inclusions $F_{v_0}' \subset F \subset F_{v_0}''$. Then the generalized
flag $\cF \cup \{F \}$ properly contains $\cF$. Conversely, if 
$\dim(F_v''/F_v') = 1$ for every nonzero $v \in V$, and if $\cG$ is a generalized 
flag which contains $\cF$, then  
$F_v' \subset G_v' \subset G_v'' \subset F_v''$. Hence $F_v' = G_v'$, $F_v'' = G_v''$, i.e.
$\cF = \cG$. 

The map $\fl$ establishes a bijection between maximal chains of subspaces
in $V$ and maximal generalized flags in $V$. More precisely, if $\cC$ is a maximal chain,
$\fl (\cC)$ is the unique maximal generalized flag which is a subchain of $\cC$.
Conversely, $\cC$ is the unique maximal chain containing $\fl (C)$. These latter statements are
essentially equivalent to Theorem 9 in \cite{DP}. For example, if $\cF$ is the maximal generalized
flag from Example 1 (iii), its corresponding maximal chain $\cC$ is described in Example 2 (iii).

We conclude this section by introducing isotropic generalized flags.
Let $w: V \times V \to V$ be a non--degenerate symmetric or skew--symmetric 
bilinear form on $V$. Denote by $U^\perp$ the 
$w$--orthogonal complement of a subspace $U \subset V$. A generalized flag $\cF$ in $V$ is
{\it $w$--isotropic} if there
exists an involution $\tau: \cF \to \cF$ such that $(F')^\perp = \tau(F)''$
for every $F' \in \cF'$. 
If $\cF$ is $w$--isotropic, $\tau$ is determined 
by $w$ and is a decreasing map. If  $\tau$ has a fixed point, the latter is
unique and will be denoted  by $F_\tau$. If $\tau$ has no fixed point,
we introduce $F_\tau' \not \in \cF$ by setting $F_\tau' := \cup_{F \subset \tau(F)} F$.
One checks immediately that $F'_\tau = \cap_{\tau(F) \subset F} F$.

\section{Compatible bases}
If $V$ is finite dimensional, any ordered basis determines a maximal flag in $V$.
Conversely, a maximal flag in $V$ determines a set of compatible bases in $V$.
More generally, if $V$ is any vector space, $\cF$ is a generalized flag in $V$ and 
$\{e_\alpha\}_{\alpha \in A}$ is a basis of $V$,
we say that $\cF$ and $\{e_\alpha\}_{\alpha \in A}$ are {\it compatible} if there exists
a strict partial order $\prec$ on $A$ (satisfying the condition stated in
the Conventions) such that 
$F_{e_\alpha}' = \Span \{e_{\beta} \, | \, \beta \prec \alpha\}$ and 
$\cF = \fl(\{F'_{e_\alpha}\}_{\alpha \in A})$.

Not every generalized flag admits a compatible basis. 
Indeed, let $V := \C[[x]]$ be the space of formal power series in the 
indeterminate $x$
and let $\cF$ denote the flag $\ldots \subset F_n \subset F_{n-1} \subset
\ldots \subset F_1 \subset F_0 = V$, where $F_n := x^n V$. Clearly, $\cF$ is a maximal flag in $V$ as 
$\dim (F_{n-1} / F_n) =1$ for all $n>0$. However, as $V$ is  uncountable dimensional, no basis
of $V$ can be compatible with the countable flag $\cF$.

The following proposition shows that the uncountability of $\dim V$ is crucial in the above example.
\begin{proposition} \label{pro2}
If $V$ is countable dimensional,
every generalized flag $\cF$ in $V$ admits a compatible basis.
\end{proposition}

\noindent
{\bf Proof.}
Assume first that $\cF$ is a maximal generalized flag in $V$. Let $\{l_i\}_{i \in \N}$ be a 
basis of $V$. Define inductively a basis $\{e_i\}_{i \in \N}$ of $V$ as follows.  Put $e_1 := l_1$.
Assuming that $e_1, \ldots, e_n$ have been constructed, choose
$e_{n+1}$ of the form $l_{n+1} + c_1 e_1 + \ldots + c_n e_n$ so that 
$F'_{e_{n+1}}$ is not among $F'_{e_1}, \ldots, F'_{e_n}$. Then, obviously, 
\begin{equation} \label{eq1}
\Span \{l_1, \ldots, l_n\} = \Span \{e_1, \ldots, e_n\}
\end{equation}
for every $n$ and the subspaces $F'_{e_n}$ are pairwise distinct. 
Furthermore, as it is not difficult to check, for every $F' \in \cF'$,
the set $F'' \backslash F'$ contains exactly one element of the basis $\{e_i\}_{i \in \N}$,
and hence $\N$ is linearly ordered in the following way: $i \prec j$ if and only if $F_{e_i}' \subset F_{e_j}'$. This
proves that $\cF$ is 
compatible with $\{e_i\}_{i \in \N}$. 

For a not necessarily maximal generalized flag $\cF$, it is enough to consider a
basis compatible with a maximal generalized flag $\cG$ containing $\cF$. Such a basis
is automatically compatible with $\cF$.
\hfill $\square$

Let $V$ be a finite or countable dimensional vector space and $w$ be a non--degenerate
symmetric or skew--symmetric bilinear form on $V$.
Define a basis of $V$ of the form $\{e_n, e^n\}$ to be {\it of type $C$} if
$w(e_i, e_j) = w(e^i, e^j) = 0$ and $w(e_i, e^j) = \delta_{i,j}$ for a 
skew--symmetric $w$. 
A basis of $V$ of the form $\{e_0 = e^0, e_n, e^n\}$ (respectively
$\{e_n, e^n\}$)  is {\it  of type $B$} (resp. {\it
of type $D$}) if
$w(e_i, e_j) = w(e^i, e^j) = 0$ and $w(e_i, e^j) = \delta_{i,j}$ for a 
symmetric $w$.
For uniformity  we will always label a basis of type $B$, $C$ or $D$  
simply as $\{e_i, e^i\}$ where we assume that in the case of $B$, $e_0 = e^0$
and $i$ runs over $\Z_+$ when $V$ is countable dimensional, or over 
a finite subset of $\Z_+$ when $V$ is finite dimensional, while in the cases of 
$C$ and $D$, $i$ runs over $\N$ or over a finite subset of $\N$. A {\it $w$--isotropic basis of $V$} is
by definition a basis of $V$ admitting an order which makes it a basis of
type $B$, $C$ or $D$. 
If $V$ is finite dimensional, $V$ admits a basis of type $B$, $C$ or $D$ if and only if
$w$ is symmetric and $V$ is odd dimensional, $w$ is skew--symmetric and then $V$ is necessarily
even dimensional, or $w$ is symmetric and $V$ is even dimensional respectively.
If $V$ is countable dimensional, by an appropriate modification of the Gram--Schmidt orthogonalization process, one shows
that if $w$ is symmetric, then $V$ admits both a basis of type $B$ and
a basis of type $D$, and if $w$ is skew--symmetric, then $V$ admits a basis of type $C$.

In the rest of the paper we assume that the dimension of $V$ is countable. 

\begin{proposition} \label{pro3}
Let $\cF$ be an $w$--isotropic generalized flag in $V$. 
\newline
{\rm(i)} Assume that $w$ is skew--symmetric. Then $V$ admits a basis of type $C$
compatible with $\cF$.
In particular, the vector space $F_\tau''/F_\tau'$ is even dimensional or infinite
dimensional.
\newline
{\rm(ii)} Assume that $k$ is algebraically closed and $w$ is symmetric.
If the vector space $F_\tau''/F_\tau'$ is odd dimensional or infinite
dimensional, then $V$ admits a  basis of type $D$ compatible with $\cF$. 
If $F_\tau''/F_\tau'$ is even dimensional or infinite
dimensional, then $V$ admits a basis of type $B$ compatible with $\cF$. 
\end{proposition}

\noindent
{\bf Proof.}
Let $\{l_n\}_{n \in \N}$ be a basis of $V$ compatible with $\cF$. 
Set $U_{F'} := \Span \{ l_i \, | \, F_{l_i}' = F'\}$ for $F' \in \cF'$.
Then $F' = \oplus_{G' \in \cF', G' \subset F', G' \neq F'} U_{G'}$ and
$F'' = F' \oplus U_{F'}$ for every $F' \in \cF'$.
It is clear therefore that the restriction of $w$ on $U_{F'} \times U_{\tau(F')}$ is 
a non--degenerate bilinear form for every $F' \in \cF'$. Furthermore, 
if $\tau$ has a fixed point, the restriction of $w$
on $U_{F_\tau'} \times U_{F_\tau'}$ is a non-degenerate skew--symmetric (respectively
symmetric) bilinear form. If $w$ is skew--symmetric, 
this implies, in particular, that $U_{F_\tau'}$, and hence $F_\tau''/F_\tau'$ is even dimensional or 
infinite dimensional. Then $U_{F_\tau'}$ admits a basis of type $C$, $B$, or $D$ depending on
whether $w$ is skew--symmetric or symmetric and on the dimension of $U_{F_\tau'}$.
Denote such a basis by $\{l'_i, {l'}^i\}$. Let, furthermore, $\{l''_i\}$ be the subset of
$\{l_n\}$ consisting of all $l_n$ with $F_{l_n}' \subset F_\tau'$. Finally,
relabel the set $\{l'_i\} \cup \{l''_i\}$ and denote the resulting set by $\{g_n\}$,
where $g_0 := l'_0 = {l'}^0$ if $\{l'_i, {l'}^i\}$ is of type $B$.

We are now ready to construct inductively the desired $w$--isotropic basis $\{e_n, e^n\}$. 
Assume that $e_i, e^i$ have been constructed for $i \leq n$.
Put $e_{n+1} := g_{n+1} - \sum_{i=1}^{n} (w(e_i, g_{n+1}) e^i + w(g_{n+1}, e^i) e_i)$.
There exists $g \in U_{\tau(F_{g_{n+1}}')}$ such that $w(e_{n+1}, g) = 1$.
Set then $e^{n+1} := g - \sum_{i=1}^{n} (w(e_i,g) e^i + w(g, e^i) e_i)$.
One checks immediately that the basis $\{e_n, e^n\}$ is $w$--isotropic and compatible
with $\cF$. 
\hfill $\square$

\section{Ind--varieties of generalized flags} \label{sec3}
For a finite dimensional $V$, two flags belong to the same connected component of the
variety of all flags in $V$ if and only if their types coincide, i.e. if the 
dimensions of the subspaces in the flags coincide. 
If $V$ is infinite dimensional the notion of type is in general not defined, and
flags, or generalized flags, can be compared using a notion of commensurability.
Such notions are well--known in the special case of subspaces of $V$, i.e.
of flags of the form $0 \subset W \subset V$, see \cite{T} and Chapter 7 of \cite{PS}.
Below we introduce a notion of commensurability for generalized flags which
in the case of subspaces reduces to a refinement of Tate's notion of
commensurability, \cite{T}.

In the rest of the paper we fix a basis $E = \{e_n\}$ of $V$. 
In the presence of a bilinear form $w$ on $V$ we fix a $w$--isotropic basis $E = \{e_n, e^n\}$, and
whenever other bases of $V$ or generalized flags in $V$ are considered they are automatically assumed to be
$w$--isotropic.
We call a generalized flag $\cF$ 
{\it weakly compatible with $E$}
if $\cF$ is compatible with a basis $L$ of $V$ such that $E \backslash (E \cap L)$ is a finite set.
Furthermore, we define two  generalized flags $\cF$ and $\cG$ in $V$
to be {\it $E$--commensurable} if both $\cF$ and $\cG$ are weakly compatible with $E$ and there exists an 
inclusion preserving bijection 
$\varphi : \cF \to \cG$ and a finite dimensional subspace $U \subset V$, such that for every $F \in \cF$
\begin{itemize}
\item[]{(i)} $F \subset  \varphi(F) + U$ and
$\varphi(F) \subset F +  U$;
\item[]{(ii)} $\dim (F \cap U) = \dim (\varphi(F) \cap U)$.
\end{itemize}
It follows immediately from the definition that any two $E$--commensurable generalized flags 
are isomorphic as ordered sets, and that two flags in a finite dimensional space
are $E$--commensurable if and only if their types coincide. (In the latter case the condition of weak compatibility with $E$ is
empty.) Furthermore,  
$E$--commensurability is an equivalence relation. Indeed, 
it is obviously reflexive and symmetric. It is also transitive.
To see this, note first that, in the definition of $E$--commensurability,
one can replace (ii) by 
\begin{itemize}
\item[]{(ii$'$)} $\dim (F/(F \cap \varphi(F))) = \dim (\varphi(F)/(F \cap \varphi(F)))$.
\end{itemize}
Consider now $\cF, \cG$ and $\cH$, such that 
$\cF$ is $E$--commensurable with $\cG$ and $\cG$ is $E$--commensurable with $\cH$. 
Let $\varphi: \cF \to \cG$ and $\psi: \cG \to \cH$ be the respective bijections, and
$U$ and $W$ be the finite dimensional subspaces of $V$ corresponding
to $\varphi$ and $\psi$ respectively. Then $\cF$ and $\cH$ satisfy (i) and 
(ii$'$) with $\psi \circ \varphi :\cF \to \cH$ and $U + W$.

\noindent
{\bf Example 3.}

\noindent
(i) Let $\cF = \{ 0 \subset F \subset V\}$ and $\cG = \{ 0 \subset G \subset V\}$.
If $F$ and $G$ are finite dimensional, then $\cF$ and $\cG$ are automatically weakly compatible with
$E$. Furthermore, $\cF$ and $\cG$ are $E$--commensurable if and only if $\dim F = \dim G$.  
If, however, $F$ and $G$ are infinite dimensional, the condition that $\cF$ and $\cG$ are 
weakly compatible with $E$ is not automatic. For example, if $F = \Span \{e_2, e_3, \ldots\}$
and $G = \{e_2 - e_1, e_3 - e_1, \ldots \}$, then $\cF$ is weakly compatible with $E$ but $\cG$ is not,
and consequently, $\cF$ and $\cG$ are not $E$--commensurable. Finally, if $F$ and $G$ are
both of finite codimension in $V$, and $\cF$ and $\cG$ are weakly $E$--compatible, then $\cF$ and $\cG$
are $E$--commensurable if and only if $\codim_V F = \codim_V G$.

\noindent
(ii) Let $\cF = \{0 = F_1 \subset F_2 \subset F_3 \subset \ldots \}$ and $\cG = \{0 = G_1 \subset G_2 
\subset G_3 \subset \ldots \}$ be two finite or infinite accending flags in $V$ compatible with $E$. 
If all subspaces $F_i$ and $G_i$
are finite dimensional, then $\cF$ and $\cG$ are $E$--commensurable if and only if
$\dim F_i = \dim G_i$ for every $i$, and $F_n = G_n$ for large enough $n$.
If, however, there are infinite dimensional spaces among $F_i$ and $G_i$, the above conditions are
still necessary for $\cF$ and $\cG$ to be $E$--commensurable but they are not always sufficient. 
The exact sufficient conditions can be derived as a consequence of the proof of Proposition 5 below.

Given a generalized flag $\cF$ weakly compatible with $E$, we denote by $\cFl(\cF, E)$
the set of all generalized flags in $V$ $E$--commensurable with $\cF$. 
For the rest of the paper we fix the following notations: 
$E_n = \{e_i\}_{i \leq n}$, $V_n = \Span E_n$, $E_n^c := \{e_i\}_{i > n}$ and $V_n^c := \Span E_n^c$. 
If $\cF$ is a $w$--isotropic generalized flag in $V$, $\cFl(\cF, w, E)$ stands for the set of all 
$w$--isotropic generalized flags $E$--commensurable with $\cF$. 
If $\cG$ is an isotropic generalized flag in $V$, the involution $\tau$ is an order reversing
involution on $\cG$ considered as an ordered set. In this case
$E_n = \{e_i, e^i\}_{i \leq n}$, $V_n = \Span E_n$, $E_n^c := \{e_i, e^i\}_{i > n}$ and $V_n^c := \Span E_n^c$. 
Since all generalized flags in $\cFl(\cF, w, E)$ are isomorphic as ordered sets, we will use the same letter
$\tau$ to denote the involution on any $\cG \in \cFl(\cF, w, E)$. 

\begin{proposition} \label{pro47}
$\cFl(\cF, E)$, as well as $\cFl(\cF, w, E)$, has a natural structure of an ind--variety. 
\end{proposition}

\noindent
{\bf Proof.} We present the proof in the case of $\cFl(\cF, E)$ only. The reader
will supply a similar proof for $\cFl(\cF, w, E)$. 
For any $\cG \in \cFl(\cF, E)$ choose a positive integer
$n_\cG$ such that $\cF$ and $\cG$ are compatible with bases containing $E_{n_\cG}^c$, 
and $V_{n_\cG}$ contains a finite dimensional subspace $U$ which
(together with the corresponding $\varphi$) makes $\cF$ and $\cG$ $E$--commensurable.
Obviously we can pick $n_\cF$ so that $n_\cF \leq n_\cG$ for every $\cG \in \cFl(\cF, E)$. 
Set also 
\begin{equation} \label{eq_new}
\cG_n := \{G \cap V_n \, | \, G \in \cG\}
\end{equation}
for $n \geq n_\cG$.

The type of the flag $\cF_n$ yields a sequence of
integers $0 = d_{n,0} < d_{n,1} < \ldots < d_{n, s_{n-1}} < d_{n, s_n} = n$,
and $\cFl(\cF_n, E_n)$ is the usual flag variety $\cFl(d_n; V_n)$
of type $d_n = (d_{n,1}, \ldots,$ $d_{n, s_n-1})$ in $V_n$.
Notice that $s_{n+1} = s_n$ or $s_{n+1} = s_n +1$. Furthermore, in both cases an integer $j_n$ is
determined as follows: in the former case, 
$d_{n+1, i} = d_{n,i}$ for $0 \leq i < j_n$ and $d_{n+1, i} = d_{n,i} + 1$ for 
$j_n \leq i < s_n$, and in the latter case $d_{n+1, i} = d_{n,i}$ for $0 \leq i < j_n$
and $d_{n+1, i} = d_{n,i-1} + 1$ for $j_n \leq i < s_n$.

Now we define a map $\iota_n : \cFl(d_n; V_n) \to \cFl(d_{n+1}, V_{n+1})$ for every $n \geq n_\cF$.
Given $\cG_n = \{0=G_0^n \subset G_1^n \subset \ldots \subset G_{s_n}^n = V_n\} \in \cFl(d_n; V_n)$, put
 $\iota_n(\cG_n) = \cG_{n+1} := \{0 = G_0^{n+1} \subset G_1^{n+1} \subset \ldots \subset G_{s_{n+1}}^{n+1}=V_{n+1}\}$,
where
\begin{equation} \label{eq2}
G_i^{n+1} := \left\{
\begin{array}{lclcl}
G_i^n & {\text { if }} & 0 \leq i < j_n && \\
G_i^n \oplus k e_{n+1} & {\text { if }} & j_n \leq i \leq s_{n+1} & {\text { and }} & s_{n+1} = s_n \\
G_{i-1}^n \oplus k e_{n+1} & {\text { if }} & j_n \leq i \leq s_{n+1} & {\text { and }} &  s_{n+1} = s_n + 1.
\end{array} \right.
\end{equation}
It is clear that 
$\iota_n$ is a closed immersion of algebraic varieties, and hence $\inj \cFl(d_n; V_n)$
is an ind--variety. Let $\psi_n : \cFl(d_n; V_n) \to \inj \cFl(d_n;V_n)$ denote the canonical
embedding corresponding to the direct system $\{ \iota_n\}$.

To endow $\cFl(\cF, E)$ with an ind--variety structure we construct
a bijection $\cFl(\cF, E) \to \inj \cFl(d_n;V_n)$.
Set
$$
\theta: \cFl(\cF, E) \to \inj \cFl(d_n;V_n), \quad \quad \theta(\cG) := \inj \cG_n,
$$
see (\ref{eq_new}).
Checking that $\theta$ is injective is straightforward.
To check that $\theta$ is surjective, fix $\widetilde{\cG} = \inj \widetilde{\cG_n} \in \inj \cFl(d_n;V_n)$,
an integer $\tilde{n}$ and a flag $\widetilde{\cG_{\tilde{n}}} \in \cFl(d_{\tilde{n}};V_{\tilde{n}})$ 
with $\psi_{\tilde{n}}(\widetilde{\cG_{\tilde{n}}}) = \widetilde{\cG}$. 
Denote by $\varphi_{\tilde{n}}$ the inclusion preserving bijection
$\varphi_{\tilde{n}}: \cF_{\tilde{n}} \to \widetilde{\cG_{\tilde{n}}}$. 
For every $F \in \cF$, put 
$\varphi(F) := \varphi_{\tilde{n}} (F \cap V_{\tilde{n}}) \oplus (F \cap V_{\tilde{n}}^c)$.
It is clear  that $\cG := \{\varphi(F)\}_{F \in \cF}$ is a generalized flag in $V$ $E$--commensurable
 with $\cF$ via $\varphi$ and $V_{\tilde{n}}$. Furthermore, using (\ref{eq2}), one verifies that
$\theta(\cG) = \widetilde{\cG}$. Hence $\theta$ is surjective.
\hfill
$\square$

\noindent
{\bf Example 4.}
 Let $\cF = \{0 \subset F \subset V\}$, see Example 3 (i). If $F$ is a finite dimensinal subspace of $V$ of
dimension $l$, then, regardless of $E$, $\cFl(\cF, E)$ is nothing but the ind--variety $Gr(l; V)$ introduced in
Section 1. If $F$ is infinite dimensional subspace of $V$ of codimension $l$, then
as a set $\cFl(\cF, E)$ depends on the choice of $E$. However, the isomorphisms
between $Gr(l; V_n)$ and $Gr(n-l; V_n)$ extend to an ind--variety isomorphism between $\cFl(\cF, E)$ and
$Gr(l; V)$ which depends on $E$. The latter isomorphism is a particular case of the following general duality.
Let $\cF$ be an arbitrary generalized flag in $V$.
Assume that $E$ is compatible with $\cF$, and for every $F \in \cF$ set $F^c := \Span\{e \in E \, | \, e \not \in F\}$.
Then $\cF^c := \{F^c \, | \, F \in \cF\}$ is a generalized flag in $V$ compatible with $E$ and
$\cFl(\cF^c, E)$ is isomorphic to $\cFl(\cF, E)$.

We complete this section by defining big cells in $\cFl(\cF, E)$ and $\cFl(\cF, w, E)$.
Let $L =\{l_n\}_{n \in \N}$ be a basis of $V$ compatible with $\cF$ and such that $E \backslash (E\cap L)$
is a finite set, and let
$U_{F'} = \Span \{l \in L \, | \, F_l' = F'\}$ for any $F' \in \cF'$.
Denote by  $\Phi = \{ \Phi_{F'}\}_{F' \in \cF'}$
a set of linear maps of finite rank $\Phi_{F'} : F' \to U_{F'}$, such that
$\Phi_{F'} \neq 0$ for finitely many subspaces $F_1' \subset \ldots \subset F_p'$ only.
Given $\Phi$, define 
$$
\begin{array}{lr}
\Gamma_{F'}: F' \to F'', & \Gamma_{F'}(v):= v + \Phi_{F'}(v), \\
\Gamma: V \to V, & \Gamma(v) := \Gamma_{F_p'} \circ \ldots \circ \Gamma_{F_{i+1}'} (v),
\end{array}
$$
where $i$ is the largest integer with $v \not \in F_i'$.
Put $\Phi(\cF) := \fl(\{\Gamma(F')\}_{F' \in \cF'})$. 
Then define the {\it big cell $C(\cF, E; L)$ of $\cFl(\cF, E)$ corresponding to the 
basis $L$} by setting 
$$
C(\cF, E; L) := \{\Phi(\cF) \, | \, {\text { for all possible }} \Phi\}.
$$ 

To define the big cell $C(\cF, w, E; L)$ in $\cFl(\cF, w, E)$, we
start with a $w$--isotropic basis $L = \{l_n, l^n\}$ of $V$ compatible with $\cF$ and
such that  $E \backslash (E\cap L)$
is a finite set,
and repeat the above construction of $\Gamma (F')$  for all $F' \in \cF'$ with
$F' \subset \tau(F')$. As a result we obtain subspaces $\Gamma(F')$ 
for $F' \subset \tau(F')$ and
set $\Phi(\cF) := \fl(\{\Gamma(F'), (\Gamma(F'))^\perp\}_{F'\in \cF, F' \subset \tau(F')})$.
Then, we set 
$$
C(\cF, w, E; L) := \{\Phi(\cF) \, | \, {\text { for all possible }} \Phi\}.
$$

Note that the role of $\cF$ in defining big cells is not special and that big cells $\cC(\cG, E; L)$,
or $\cC(\cG, w, E; L)$, are well defined for every $\cG \in \cFl (\cF, E)$, or respectively
$\cG \in \cFl (\cF, w, E)$.

\begin{proposition} \label{pro45} $\phantom{xxxx}$

\noindent
{\rm (i)} The big cell $C(\cF, E; L)$ (respectively, $C(\cF, w, E; L)$)
is an affine open  ind--subvariety of $\cFl(\cF, E)$ (resp. of $\cFl(\cF, w, E)$).

\noindent
{\rm (ii)} We have
\begin{equation} \label{eqA1}
\cFl(\cF, E) = \cup_L \,\, C(\cF, E; L)
\end{equation}
and
\begin{equation} \label{eqA2}
\cFl(\cF, w, E) = \cup_L \,\, C(\cF, w, E;  L),
\end{equation}
where the unions run over all bases (respectively $w$--isotropic bases) $L$ of $V$ compatible 
with $\cF$ and such that $E \backslash (E \cap L)$ is a finite set.
\end{proposition}

\noindent
{\bf Proof.} 
We discuss the case of $\cFl(\cF, E)$ only. The argument for the 
case of $\cFl(\cF, w, E)$ is similar.
Put $L_n := \{l_i\}_{i \leq n}$ and $W_n := \Span L_n$.
Let  $\cF_n$ and $\cFl(d_n; W_n)$ be as in the proof of 
Proposition \ref{pro47} above.
Set $\cC(d_n; W_n; L_n) := \{\Phi(\cF)_n \, | \, {\text {for all }} \Phi {\text { such that, for every }} 
F' \in \cF',  \Phi_{F'}(W_n)$  $\subset W_n 
{\text { and }} \Phi_{F'}(l_i)$ $ =0 {\text { for }} i > n \}$.
Obviously, $\cC(d_n; W_n; L_n)$ is a big cell in $\cFl(d_n; W_n)$, and 
hence is an affine open subset. 
Therefore, the inclusion $\iota_n (\cC(d_n; W_n; L_n))$ $\subset 
\cC(d_{n+1}; W_{n+1}; L_{n+1})$ and the equality 
$\inj \cC(d_n; W_n; L_n) = C(\cF, E; L)$ show that
$C(\cF, E; L)$ is an affine open ind--subvariety of $\cFl(\cF, E)$. 
The fact that the set of cells $\{C(\cF, E; L) \, | \, L$ is a basis of $V 
{\text { compatible with }} \cF, {\text { such that }} L \backslash (E \cap L) {\text { is a finite set}}\}$
is a covering of $\cFl(\cF, E)$
is an easy consequence of the definition of $E$--commensurability.
\hfill $\square$
 
\section{Ind--varieties of generalized flags as homogeneous ind--spaces} \label{sec4}
Let $G(E)$ be the group of automorphisms $g$ of $V$ 
such that $g(e) = e$ for all but finitely many $e \in E$ and in addition $\det g = 1$.
Recall that $E_n = \{e_i\}_{i \leq n}$ and
$V_n = \Span E_n$. The natural inclusion 
$$
G(E_n) \subset G(E_{n+1}), \quad \quad g \mapsto \kappa_n (g),
$$
where $\kappa_n(g)_{|V_n} = g$ and $\kappa_n(g) (e) = e$ for $e \in E_{n+1} \backslash E_n$,
is a closed immersion of algebraic groups. Furthermore, $G(E) = \cup_{n \in \N} G(E_n)$.
In particular $G(E)$ is a locally linear ind--group, and $G(E) = G(L)$ for any basis $L$ of $V$ such that $E \backslash (E \cap L)$
is a finite set.

Similarly, when $E$ is a $w$--isotropic basis of $V$, let 
$G^w(E) := \{ g \in G(E) \,|\, w(g(u), g(v)) = w(u, v) {\text { for any }} u, v \in V\}$. 
There are natural closed immersions $G^w(E_n) \subset G^w(E_{n+1})$, 
and $G^w(E) = \cup_{n \in \N} G^w(E_n)$, where in this case $E_n := \{e_i, e^i \}_{i \leq n}$.

The ind--group $G(E)$ (respectively, $G^w(E)$) is immediately seen to be isomorphic to
the classical ind--group $A(\infty)$ (resp., $B(\infty)$, $C(\infty)$ or $D(\infty)$ 
if $E$ is a $w$--isotropic basis of type $B$, $C$ or $D$).
The ind--groups $A(\infty)$, $B(\infty)$, $C(\infty)$ and $D(\infty)$ are discussed in
detail in \cite{DPW}. An alternative notation for $A(\infty)$ is $SL(\infty)$, and 
$B(\infty) \cong D(\infty)$ and $C(\infty)$ are also denoted respectively by $SO(\infty)$
and $Sp(\infty)$.

In the rest of the paper the letter $G$ will denote one of the groups $G(E)$ or 
$G^w(E)$, and $G_n$ will denote 
respectively $G(E_n)$ or $G^w(E_n)$.
The basis $E$ equips $G$ with a subgroup $H$, consisting of all diagonal automorphisms of
$V$ in $G$, i.e. of the elements $g \in G$ such that $g(e) \in k e$ for every $e \in E$.
We call $H$ a {\it splitting Cartan subgroup} 
(in the terminology of \cite{DPW}, $H$ is a Cartan subgroup of $G$). 
Following \cite{DPW}, for the purposes of the present paper, we define a 
{\it parabolic} (respectively, a {\it Borel}) {\it subgroup of $G$} to be
an ind--subgroup $P$ (resp., $B$) of $G$  such that its intersection with $G_n$ 
 for every $n$ is a parabolic (resp., a Borel) subgroup of $G_n$ for some,
or equivalently any, order on $E$. 

If $\cF$ is a generalized flag in $V$ compatible with $E$ (and $w$--isotropic, 
whenever $E$ is $w$--isotropic), we denote by $P_\cF$ the stabilizer of $\cF$
in $G$.

\begin{proposition} \label{pro5} $\phantom{xxxx}$ \newline
{\rm (i)} $P_\cF$ is a parabolic subgroup of $G$ containing $H$.
\newline
{\rm (ii)} The map $\cF \mapsto P_\cF$ establishes a bijection between generalized flags
in $V$ compatible with $E$ and parabolic subgroups of $G$ containing $H$. 
\end{proposition}

\noindent
{\bf Proof.} The inclusion $H \subset P_\cF$ follows directly 
from the definition of $H$ and $P_\cF$. 
Furthermore, $P_\cF \cap G_n$ is a parabolic subgroup of $G_n$ as it is 
the stabilizer of $\cF_n$ in $G_n$. Hence $P_\cF$ is a parabolic subgroup of $G$.
If, conversely, $P = \cup_n P_n$ is a parabolic subgroup of $G$ containing $H$,
denote by $\cF(n)$ the flag in $V_n$ whose stabilizer is $P_n$. 
Note that $\cF(n)$ maps into $\cF(n+1)$. More precisely, for $G = G(E)$,
$\cF(n+1) = \iota_n (\cF(n))$, see (\ref{eq2}); 
and for $G = G^w(E)$, the corresponding map is the $w$--isotropic analog
of $\iota_n$ which we leave to the reader to reconstruct. 
In both cases we define $\cF$ as $\theta^{-1} (\inj \cF(n))$.
A direct checking shows that $P = P_\cF$.
\hfill $\square$

Proposition \ref{pro5} further justifies our 
consideration of generalized flags, see the discussion before Proposition \ref{pro02} above. 
Indeed, it is clear that if $P \subset G$ is
the stabilizer of a chain $\cC$ of subspaces in $V$, then $P$ depends only
the partition $\pi(\cC)$, see (\ref{eq01}), and not on $\cC$ itself. Therefore,
the generalized flag $\cF$ emerges as a representative of the class
of all chains $\cC$ which have $P$ as a stabilizer in $G$. Moreover,
Proposition \ref{pro2} together with Proposition \ref{pro5} 
(respectively, Proposition \ref{pro3} and Proposition \ref{pro5} for 
$w$--isotropic flags), 
imply that the stabilizer in $G$ of any generalized flag (resp. isotropic generalized flag)
compatible with $E$  is a parabolic subgroup of $G$.
Finally,  maximal generalized flags in $V$ correspond to Borel subgroups
under the above bijection.

Note that, for any order on $E$ and for any generalized flag $\cF$ compatible with $E$,
$G/P_\cF = \cup_n (G_n/P_n)$, where $P_n := P_\cF \cap G_n$. In particular, $G/ P_\cF$ 
is an ind--variety.  Moreover, any other order on $E$, for which $E$ is isomorphic to $\N$
as an ordered set, defines an isomorphic ind--variety.
We are now ready to exhibit the  homogeneous ind--space structure on
$\cFl(\cF, E)$ and $\cFl(\cF, w, E)$. 

\begin{theorem} \label{the1} For any $E$ and $\cF$ as above there is a respective isomorphism
of ind--varieties $\cFl (\cF, E)  \cong G/P_\cF$ or $\cFl (\cF, w, E)  \cong G/P_\cF$.
\end{theorem}

\noindent
{\bf Proof.} 
Given $\cG \in \cFl (\cF, E)$ (or, respectively,  $\cG \in \cFl (\cF, w, E)$), let $U \subset V$ be 
the corresponding to $\cG$ finite
dimensional subspace. We may assume that $U = V_n = \Span E_n$ for some $n$.  
Since $\cF_n$ and $\cG_n$ are 
flags of the same type in the finite dimensional space $V_n$, there exists  
$g_n \in G_n$, so that $g(\cF_n) = \cG_n$. 
We extend $g_n$ to an element $g \in G$ by setting $g(e) = e$ for $e \in E \backslash E_n$.
Now
$$
f: \cFl (\cF, E) \to G/P_\cF {\text { (or, respectively, }} 
f: \cFl (\cF, w, E) \to G/P_\cF {\text {)}}, \quad f(\cG):= gP
$$
is a well--defined map and it is easy to check that it is an isomorphism of ind--varieties.
\hfill $\square$

\section{Picard group and projectivity} \label{sec5}
The interpretation of $\cFl(\cF, E)$ and $\cFl(\cF, w, E)$ as 
homogeneous ind--spaces $G/P_\cF$ provides us with a representation theoretic
description of the Picard groups of $\cFl(\cF, E)$ and $\cFl(\cF, w, E)$.
Namely, $\Pic \cFl(\cF, E)$, as well as $\Pic \cFl(\cF, w, E)$, is  naturally isomorphic to
the group of integral characters of the Lie algebra of the ind--group $P_\cF$.

Consider $\cFl(\cF, E)$. There is a canonical isomorphism of abelian groups 
$\Pic \cFl(\cF, E) = \Hom (P_\cF, k^\times)$. 
To see this, notice that $\Pic \cFl(\cF, E) =
\proj \Pic \cFl(d_n;V_n) = \proj \Pic G_n/ (P_\cF)_n$. It is a classical fact that 
$\Pic G_n/ (P_\cF)_n = \Hom((P_\cF)_n, k^\times)$ for every $n$, 
and an immediate verification shows that the diagram 
\begin{equation} \label{hhh}
\begin{array}{ccc}
\Pic (G_{n+1}/ (P_\cF)_{n+1})  & \cong & \Hom((P_\cF)_{n+1}, k^\times)\\
\downarrow                       &       & \downarrow\\
    \Pic (G_n/ (P_\cF)_n)     & \cong & \Hom((P_\cF)_n, k^\times)
\end{array}
\end{equation}
is commutative. Hence $\Pic \cFl(\cF,  E) \cong \Hom(P, k^\times)$, and $\Hom(P_\cF, k^\times)$
is nothing but the group of integral characters of the Lie algebra of $P_\cF$.
In the case of $\cFl(\cF, w, E)$ the desired isomorphism is established by replacing 
$\Hom((P_\cF)_n, k^\times)$ and $\Hom((P_\cF)_{n+1}, k^\times)$ in diagram (\ref{hhh}) 
with the groups of integral characters of the Lie algebras of 
$P_{\cF_n}$ and $P_{\cF_{n+1}}$ respectively. 

In the rest of this section we give a purely geometric description of 
$\Pic \cFl(\cF, E)$ and $\Pic \cFl(\cF, w, E)$.
Consider the corresponding covering (\ref{eqA1}) or (\ref{eqA2}).
Let $L$ and $M$ be two bases compatible with $\cF$ for which
$E \backslash (E \cap L)$ and $E \backslash (E \cap M)$ are finite sets. 
Denote by $g_{L, M}$ the automorphism of $V$  such that $g_{L, M}(l_i) = m_i$ 
for $\cFl (\cF, E)$, and $g_{L, M}(l_i) = m_i$, $g_{L, M}(l^i) = m^i$ for
  $\cFl(\cF, w, E)$. 
It has a well--defined determinant and, moreover, it induces an automorphism of
$F''/F'$. Denote the determinant of this latter automorphism by 
$\det_{L, M} (F''/F')$. In this way we obtain an invertible sheaf
$\cL_{F'}$ with transition functions $\det_{L, M}(F''/F')$ on 
$C(\cF, E; L) \cap C(\cF, E; M)$ or $C(\cF, w, E; L) \cap C(\cF, w, E; M)$
respectively. Finally, let 
$\gamma_{F'} \in \Pic \cFl(d_n;V_n)$, respectively 
$\gamma_{F'} \in \Pic$  $\cFl(d_n, w; V_n)$, denote the class of $\cL_{F'}$.

\begin{proposition} \label{pro7} There are canonical isomorphisms of abelian groups
$\Pic \cFl(\cF, E)$ 
$\cong$ \newline $(\Pi_{F' \in \cF'} (\Z   \gamma_{F'})) / (\Z  \, \Pi_{F' \in \cF'} \, \gamma_{F'})$
and
$\Pic \cFl(\cF, w, E) 
\cong \Pi_{F' \in \cF', F' \subset \tau(F'), F' \neq F_\tau'} \, (\Z  \gamma_{F'})$.
\end{proposition}

\noindent
{\bf Proof.} Consider the case of $\cFl(\cF, E)$ first.
Let $\gamma_{F',n}$ the class of the restriction $(\cL_{F'})_n$ of $\cL_{F'}$ 
to $\cFl(d_n; V_n)$. Then $\gamma_{F',n} = 0$ unless $F'' \cap V_n \neq F' \cap V_n$.
Define the group homomorphism 
$\varphi_n : \Pi_{F' \in \cF'} (\Z   \gamma_{F'}) \to \Pic \cFl (d_n; V_n)$
via $\varphi_n(\Pi_{F' \in \cF'} m_{F'}  \gamma_{F'}) := \sum_{F' \in \cF'} m_{F'} \gamma_{F', n}$.
The sum   $\sum_{F' \in \cF'} m_{F'} \gamma_{F', n}$ makes sense because $\gamma_{F', n} = 0$ for all but finitely
many $F' \in \cF'$.
Clearly $\varphi_n = r_n \circ \varphi_{n+1}$, where
$r_n: \Pic \cFl (d_{n+1}; V_{n+1}) \to \Pic \cFl (d_n; V_n)$
is the restriction map. Therefore, 
 by the universality property of $\proj$, there is a
homomorphism 
$$
\varphi: \Pi_{F' \in \cF'} (\Z   \gamma_{F'}) \to \Pic \cFl (\cF, E)
= \proj \Pic \cFl (d_n; V_n).
$$
 Furthermore $\varphi$ is surjective as $\varphi_n$ is surjective 
for each $n$. 

To compute $\ker \varphi$, note that $\ker \varphi = \cap \ker \varphi_n$. 
We have $\ker \varphi_n = (\Z (\Pi_{F' \in \cF'} \gamma_{F'})) 
\times \Pi_{F'\in \cF', F' \cap V_n = F'' \cap V_n} (\Z \gamma_{F'})$
and therefore $\ker \varphi = \Z (\Pi_{F' \in \cF'} \gamma_{F'})$,
i.e. $\Pic \cFl(\cF, E)$ 
$\cong$ \newline 
$(\Pi_{F' \in \cF'} (\Z   \gamma_{F'})) / (\Z  \, \Pi_{F' \in \cF'} \, \gamma_{F'})$.

In the case of $\cFl(\cF, w, E)$ homomorphisms 
$\varphi_n : \Pi_{F' \in \cF'} (\Z   \gamma_{F'}) \to \Pic \cFl (d_n, w; V_n)$
and $\varphi: \Pi_{F' \in \cF'} (\Z   \gamma_{F'}) \to \Pic \cFl (\cF, w, E)$ 
are defined in a similar way.
Here  $\ker \varphi_n$ $=$ $\Pi_{F' \in \cF', F' \subset \tau(F')}$ 
$(\Z(\gamma_{F'} + \gamma_{\tau(F')'}))  
\times \Pi_{F'\in \cF', F' \subset \tau(F'), F' \cap V_n = F'' \cap V_n} 
(\Z \gamma_{F'})$, and consequently
$\ker \varphi$ $=$ $\Pi_{F' \in \cF', F' \subset \tau(F')}$ $(\Z(\gamma_{F'} + \gamma_{\tau(F')'}))$ 
i.e. $\Pic \cFl(\cF, w, E) 
\cong \Pi_{F' \in \cF', F' \subset \tau(F'), F' \neq F_\tau'} \, (\Z  \gamma_{F'})$.
\hfill $\square$

We complete this paper by an explicit criterion for the projectivity of $\cFl(\cF, E)$ 
and $\cFl(\cF, w, E)$. The following Proposition is a translation of Proposition 15.1
in \cite{DPW} into the language of generalized flags.

\begin{proposition} \label{pro8} $\cFl(\cF, E)$ or 
$\cFl(\cF, w, E)$ is projective if and only if $\cF$ is a flag.
\end{proposition}

\noindent
{\bf Proof.} Consider the case of $\cFl(\cF, E)$ (the case of $\cFl(\cF, w, E)$ is
similar). 
$\cFl(\cF, E)$ is projective if and only if it admits a very ample invertible
sheaf. An immediate verification shows that an invertible sheaf $\cL$, whose class in 
$\Pic \cFl(\cF, E)$ is the image of $\Pi_{F' \in \cF'} m_{F'} \gamma_{F'}$, is
very ample if and only if the map $c: \cF' \to \Z$, $F' \mapsto m_{F'}$ is 
strictly increasing. 
Indeed, $\cFl(\cF, E) = \inj \cFl(d_n;V_n)$, 
and $\cL$ is very ample if and only if its restrictions $\cL_n$ onto $\cFl (d_n;V_n)$
are very ample for all $n$. Consider the map
$c_n : \cF_n' \to \Z$, defined via
$c_n ((\cF_n)_v'):= c((\psi_n(\cF_n))_v')$ for every nonzero $v \in V_n$. (As the reader will check,
$c_n$ is well defined, i.e. if $(\cF_n)_{v_1}' = (\cF_n)_{v_2}'$, then $(\psi_n(\cF_n))_{v_1}'
= (\psi_n(\cF_n))_{v_2}'$.)
According to the classical Bott--Borel--Weil Theorem,
$\cL_n$ is very ample if and only if the map $c_n$
is strictly increasing.    
Hence, $\cL$ is very ample if and only if $c$ is strictly increasing.
This enables us to conclude that 
$\cFl(\cF,  E)$ is projective if and only if there exists a strictly
increasing map $\cF' \to  \Z$, i.e. if and only if $\cF$ is a flag.
\hfill $\square$

Propositions 8 and 9 allow us to make some initial remarks concerning the isomorphism classes of the ind--varieties
$\cFl(\cF, E)$ and $\cFl(\cF, w, E)$. For example, if $\cF$ is a flag of finite length in $V$, and
$\cG$ is a flag (or generalized flag) in $V$ of length different from the length of $\cF$ (finite or infinite), 
then $\cFl(\cF, E)$ and $\cFl(\cG, L)$ are not isomorphic because their Picard groups are not isomorphic.
Furthermore, if $\cF$ is a flag in $V$ but $\cG$ is not, 
then $\cFl(\cF, E)$ and $\cFl(\cG, L)$ are not isomorphic because the former ind--variety is projective and the
latter is not. Finally,  a recent result of J. Donin and the second named author, \cite{DoP},
implies that if $\cF = \{0 \subset F \subset V\}$ with $F$ both infinite dimensional and of
infinite codimension in $V$, then the ``ind--grassmannian'' $\cFl(\cF, E)$ is not isomorphic to $Gr(l; V)$ for any $l$, cf. Example 4.

\vskip 2 cm

\begin{tabular}{lll}
I.D.: & \quad \quad \quad \quad \quad & I.P.: \\
Department of Mathematics and Statistics && Department of Mathematics  \\
Queen's University  && University of California   \\
Kingston, Ontario, K7L 3N6      && Riverside, CA 92521 \\
Canada && USA\\
{\tt dimitrov@mast.queensu.ca} && {\tt penkov@math.ucr.edu} 
\end{tabular}


\begin{thebibliography}{DPW}
\bibitem[DP]{DP}
I. Dimitrov and I. Penkov, Weight modules of direct limit Lie algebras,
Internat. Math. Res. Notices, No 5, 1999, 223--249.

\bibitem[DPW]{DPW}
I. Dimitrov, I. Penkov and J. A. Wolf, A Bott--Borel--Weil theory 
for direct limits of algebraic groups, Amer. J. of Math. {\bf 124} (2002), 955--998.

\bibitem[DoP]{DoP}
J. Donin and I. Penkov, Finite rank vector bundles on inductive limits of Grassmannians, 
Internat. Math. Res. Notices, 2003, No 34, 1871--1887.

\bibitem[NRW]{NRW}
L. Natarajan, E. Rodr\'iguez-Carrington and J. A. Wolf, The Bott--Borel--Weil theorem for direct limit groups,
Trans. Amer. Math. Soc. 353 (2001), 4583--4622.


\bibitem[O1]{O1} 
G. I. Olshanskii , Unitary representations of the infinite--dimensional classical groups 
${\rm U}(p,\,\infty )$, ${\rm SO}\sb{0}(p,\,\infty )$, ${\rm
Sp}(p,\,\infty )$, and of the corresponding motion groups, (Russian) Funktsional. Anal. i Prilozhen. 12 (1978), 
no. 3, 32--44, 96. 

\bibitem[O2]{O2}
G. I. Olshanskii, Description of unitary representations with highest 
weight for the groups ${\rm U}(p,\,q)\,\sp{\sim }$, (Russian) Funktsional. Anal. i
Prilozhen. 14 (1980), no. 3, 32--44, 96.


\bibitem[PS]{PS}
A. Presley and G. Segal, Loop groups, Clarendon Press, Oxford, 1986.

\bibitem[SV]{SV}
\c S. Str\u atil\u a and D.  Voiculescu,  A survey on representations of the unitary group ${\rm U}(\infty )$,
Spectral theory (Warsaw, 1977), 415--434, Banach Center Publ., 8, PWN, Warsaw, 1982.

\bibitem[T]{T}
J. Tate, Residues of differentials on curves, Ann. csient. \'Ec. Norm. Sup.,
4$^e$ s\'erie, {\bf 1} (1968), 149--159.

\bibitem[VK1]{VK1}
A. M. Vershik and S. V.  Kerov, Characters and factor-representations of the infinite unitary group, 
(Russian) Dokl. Akad. Nauk SSSR 267 (1982), no. 2, 272--276.

\bibitem[VK2]{VK2}
A. M. Vershik and S. V. Kerov,  Locally semisimple algebras. Combinatorial theory and the 
$K\sb 0$-functor, (Russian) Current problems in mathematics. Newest
results, Vol. 26, 3--56, 260, Itogi Nauki i Tekhniki, Akad. Nauk SSSR, 
Vsesoyuz. Inst. Nauchn. i Tekhn. Inform., Moscow, 1985.





\end{thebibliography}
\end{document}